\documentclass[10pt,twoside]{article}

\usepackage{amsfonts} 
\usepackage{amsmath, color}
\usepackage{pstricks}  

\usepackage{graphicx}
\usepackage{esint}

\usepackage{amssymb}
\usepackage{latexsym}
\usepackage{hyperref}

\makeatletter
\newtheorem{theorem}{Theorem}[section]
\newtheorem{remark}[theorem]{Remark}
\newtheorem{lemma}[theorem]{Lemma}
\newtheorem{proposition}[theorem]{Proposition}
\newtheorem{corollary}[theorem]{Corollary}
\newtheorem{definition}[theorem]{Definition}
\newtheorem{example}[theorem]{Example}

\def\Imm{\rm Imm}
\def\1{\mathbf{1}}
\def\:{\lrcorner}
\def\#{\sharp}

\def\a{\alpha}
\def\b{\beta}
\def\g{\gamma}

\def\e{\varepsilon}

\def\o{\circ}
\def\s{\sigma}

\def\x{\otimes}

\def\qed{\ensuremath{\quad\blacksquare\quad}}

\def\inv#1{\raise.1em\hbox to 0pt{$^{-1}$\hss}_{#1}\;}

\def\V{\noindent}
\def\v{\noindent}

\newcommand{\bean}{\begin{eqnarray*}}
\newcommand{\eean}{\end{eqnarray*}}
\newcommand{\benu}{\begin{enumerate}}
\newcommand{\eenu}{\end{enumerate}}
\newcommand{\eea}{\end{eqnarray}}
\newcommand{\bea}{\begin{eqnarray}}

\newcommand{\coker}{{\rm coker}}

\newcommand{\tr}{{\rm tr}}

\def \beit{\begin{itemize}}
\def \eeit{\end{itemize}}
\def \bey{\begin{eqnarray*}}
\def \eey{\end{eqnarray*}}

\def \bui#1#2{\mathrel{\mathop{\kern 0pt#1}\limits^{#2}}}
\def \buil#1#2{\mathrel{\mathop{\kern 0pt#1}\limits_{#2}}}

\newtheorem{Theorem}{Theorem}

\newtheorem{Lemma}{Lemma}

\newcommand{\be}{\begin{equation}}
\newcommand{\ee}{\end{equation}}

\newcommand{\N}{{\mathbb N}}
\newcommand{\Z}{{\mathbb Z}}
\newcommand{\R}{{\mathbb R}}

\newcommand{\G}{\Gamma}

\newcommand{\ben}{\begin{enumerate}}
\newcommand{\een}{\end{enumerate}}
\newcommand{\bit}{\begin{itemize}}
\newcommand{\eit}{\end{itemize}}
\newcommand{\edoc}{\end{document}}

\newcommand{\bdefi}{\begin{definition}}
\newcommand{\btheo}{\begin{theorem}}
\newcommand{\bprop}{\begin{proposition}}
\newcommand{\brema}{\begin{remark}}
\newcommand{\bcoro}{\begin{corollary}}
\newcommand{\blemm}{\begin{lemma}}
\newcommand{\bexam}{\begin{example}}

\newcommand{\edefi}{\end{definition}}
\newcommand{\etheo}{\end{theorem}}
\newcommand{\eprop}{\end{proposition}}
\newcommand{\erema}{\end{remark}}
\newcommand{\ecoro}{\end{corollary}}
\newcommand{\elemm}{\end{lemma}}
\newcommand{\eexam}{\end{example}}

\title{Applying the index theorem to non-smooth operators}

\begin{document}

\author{Olaf M\"uller\footnote{Fakult\"at f\"ur Mathematik, Universit\"at Regensburg, D-93040 Regensburg, \texttt{Email: olaf.mueller@ur.de}}}

\date{\today}
\maketitle

\begin{abstract}
\v We give a simple way to extend index-theoretical statements from partial differential operators with smooth coefficients to operators with coefficients of finite Sobolev order. \end{abstract}

In the theory of partial differential equations and in particular in index theory it is most common to consider only partial differential operators $P$ with smooth coefficients, or, more generally, smooth pseudo-differential operators, which still map smooth sections to smooth sections. If $P$ is moreover elliptic, this allows defining a parametrix of $P$, that is, an inverse up to infinitely smoothing operators. However, sometimes the assumption of smooth coefficients is not justified, e.g. in the theory of geodesic curves on some natural spaces of mappings. Here, to ensure some minimal amount of computibility, one usually restricts the considerations to spaces of mappings of finite Sobolev or H\"older regularity instead of smooth mappings. Consequently, the induced metrics are also of finite regularity only, and so are the differential operators defined in terms of them, which appear in the geodesic equations. Often, as in \cite{BHM}, the proof of local existence and uniqueness for the geodesic equations requires an application of the index theorem to such an operator with Sobolev coefficients seen as as a Fredholm map between Sobolev spaces $H^{k+p} (\pi)  \rightarrow H^k (\pi) $ of sections of a vector bundle $\pi: E \rightarrow M$ over a compact manifold $M$. One could hope that the notion of abstract elliptic operators due to Atiyah \cite{mA} and elaborated on topological manifolds by Teleman \cite{nT} embark the case of partial differential operators with Sobolev coefficients. Amazingly enough, this does not seem to be the case, as one of the assumptions of the theorems in those articles is that multiplication with every $f \in C^0(M) $ is a bounded operator on the Hilbert spaces used, which is not true for any Sobolev Hilbert space except $L^2(M)$. Here we want to present a solution to the above problem, assuming coefficients of regularity increasing linearly with the order, a condition satisfied by the Laplacian. After the prepublication of the first version of the article, the author learned about the techniques presented in the third chapter of \cite{mT2}. It would be interesting to compare the results of this note to the results about special first and second order pseudodifferential elliptic operators obtained there. The author wants to thank Ulrich Bunke for an interesting discussion on the topic of this note and Martin Bauer, Martins Bruveris, Peter Michor and an anonymous referee for many useful suggestions.

{\bf Conventions:} We impose the convention $0 \in \N$, and for $T \in \N$ we define $\N_T:= \N \cap [1,T]$. Throughout the article, $M$ denotes an arbitrary smooth compact $n$-dimensional manifold. Contrary e.g. to \cite{mS}, no smoothness assumption is made on the coefficients of a differential operator (as functions on the corresponding jet bundles) --- they are a priori only assumed to be measurable.   

\newpage

Let $M$ be a smooth compact manifold of dimension $n$, let $k, N \in \N$, and let $\pi$ be an $H^k$ vector bundle of rank $N$ over $M$, i.e. there is a bundle atlas $(U_t, D_t)$, $t \in \N_T$, with $H^k$ trivialization changes. If $g$ is an $H^k$ vector bundle scalar product on $\pi$, $\omega $ is a volume form on $M$ and $\nabla$ is a metric connection on $\pi$ with Christoffel symbols of regularity $H^{k-1}$ in the bundle charts, then $(\pi, g, \nabla)$ is called a {\bf metric $H^k$ vector bundle with connection}. In this case we define in the usual way the norm $W^{l,q}(g, \nabla)$ (and the scalar product $H^l (g, \nabla) $ associated to $W^{l,2}(g, \nabla)$) on the Sobolev space $W^{l,q}(\pi)$ of $W^{l,q}$ sections of $\pi$ for any $l \leq k$. A differential operator $P$ of order $s \leq k$ on $\pi$ is called {\bf $k$-safe} iff for every $U_t$ and every multiindex $i = (i_1, ... , i_r)$ there are $(N \times N)$-matrices $P_{i,t} $ with 

\bean
({\rm jet}^s D_t) \o P \o ( D_t)^{-1} = \sum_{r=0}^s \sum_{1 \leq ... i_1 \leq i_2 \leq ... \leq i_r \leq n} P_{i,t} \partial_{i_1} ... \partial_{i_r},
\eean

and such that every matrix $P_{i,t}$ is of regularity $H^{p(i,k)}$ where 
\bean
p(i,k) := {\rm max} \{ k - s, \vert i \vert-s + \lfloor \frac{n}{2} \rfloor + 1 \} ,
\eean

\v where, for $i \in \N$, $\lfloor i \rfloor$ is the Gauss bracket of $i$. The first thing to note is that $ p(i,k) >\frac{n}{2} $ for all $i$ with $\vert i \vert = s$, i.e., the leading coefficients are continuous, which means that ellipticity can be defined as usual. Obviously, the $k$-safe operators of degree $s$ form a linear subspace $S_k^s (\pi)$ of the space of operators of degree $s$, which we equip with the scalar product $(\cdot , \cdot ) $ defined by 

\bean
(P,Q):= \sum_{t \in \N_T }   \sum_{\vert i \vert \leq s}  H^{a(i,k)} ({\rm pr}_2^* h, \nabla_0)(P_{i,t} , Q_{i,t} ) 
\eean

where $h $ is any scalar product on $\R^{N^2}$ (any two choices yield topologically equivalent scalar products), $\nabla_0$ is the trivial connection. There are more invariant ways to define a scalar product of this kind, but the above definition will fully suffice for our purposes. We will often write $S^s_k$ instead of $S_k^s (\pi)$. By definition, we have $S_{k+1}^s \subset S_k^s$. Obviously, if $\pi_1$ is a metric $H^{k(1)}$-vector bundle with connection and $\pi_2 $ is a metric $H^{k(2)}$-vector bundle with connection then $\pi_1 \oplus \pi_2$ and $\pi_1 \x \pi_2$ with the corresponding natural metrics and connections are metric $H^k$-vector bundles with connections, for $k:= \min \{ k_1, k_2 \} $. To make an analogous statement on the dual vector bundle, recall that on a compact manifold standard estimates (see e.g. \cite{CD}, VI.3 (p. 386)) give $W^{k,p} \cdot W^{l,q} \subset W^{m,p}$ for all $p, q \in \R$ with $  1 < p \leq q \leq \infty $ and  for all $k,l,m \in \N$ with $ k+l > m + n/q $ and $  k,l \geq m $. Specializing to Sobolev Hilbert spaces ($p=q=2$), this means $H^k \cdot H^l \subset H^{S(k,l)} $ for 

\bea
\label{RingPropertyHilbert}
S(k,l) := \min \big\{ k,l, k+l - \lfloor n/2 \rfloor -1  \big\} ,
\eea

as long as $k+l - \lfloor n/2 \rfloor -1 \geq 0 $.
\begin{Theorem}
\label{CStarAlgebra}
Let $P$ be a differential operator of order $s$ on a vector bundle $\pi:E \rightarrow M$, let $l > n/2$.
\begin{enumerate}
\item $P^{(l)}:= P \vert_{H^l (\pi)}: H^l (\pi) \rightarrow H^{l-s} (\pi)$ is bounded if and only if $P$ is $l$-safe, and $e^{(l)} : P \mapsto P^{(l)}$ is a smooth linear map from $S^s_l (\pi)$ to $(H^l)' (\pi)) \x H^{l-s} (\pi)$.
\item  If $k>n/2$ and if $\pi$ is equipped with a vector bundle metric $g$ of regularity $H^k$ and $M$ has a volume form $\omega$ of regularity $H^{k-s}$, then for $k-s >n/2$, taking the formal $H^0(g, \omega)$-adjoint $P'_{g, \omega}$ of a $k$-safe operator $P$ of order $s$ is a bounded map $S_k^s (\pi)\rightarrow S_{k-s}^s (\pi)$.
\item {If $k  > n/2$ and $s_1+ s_2 \leq k$, then composition of operators is a bounded bilinear map 

$S_k^{s_1}(\pi) \times S_k^{s_2} (\pi)\rightarrow S_k^{s_1 + s_2} (\pi)$.}
\end{enumerate}
\end{Theorem}

{\bf Proof.}  The first assertion can be verified by simply counting orders. The necessity of the condition can be seen by the classical fact of sharpness of $S(u,v)$ in Eq. (\ref{RingPropertyHilbert}) (which can be verified in an open ball in $\R^n$ considering the functions $x \mapsto \vert x \vert^\a$, for different values of $\a$) combined with appropriate partitions of unity.

Now the formal $H^0(g)$-adjoint $P'_{g, \omega}$ of an operator $P:= \sum_{\vert j \vert \leq s} P_j \partial_{j}$ (by linearity w.l.o.g. applied to a section $u$ supported in a single trivialization chart) is given as 

$$ (P'_{g,\omega}) u(x) = \sum_{\vert j \vert <s} (-1)^{\vert j \vert} \partial_j (P^g_j (x) w(x) u(x)) ,$$

where $P^g$ denotes the transpose of a matrix with respect to $g$ and $\omega = w dx_1 \wedge ... \wedge dx_n$ in the local coordinates. Thus for the coefficients of $P'_{g, \omega}$ we get 

$$ (P'_{g, \omega})_{j} = \sum_{ l \geq j} (-1)^{\vert j \vert} \binom{l}{j} \partial_{l-j} (P^g_l \cdot w),$$

where we used subtraction of multiindices. If $P$ is $k$-safe, then $P_l, P_l^g \in H^{a(l,k) }$, and $P'_{g, \omega}$ is seen to be $(k-s)$-safe by counting orders. It is easy to see that the map is continuous. 

For the third assertion, one can invoke the equivalence of the first assertion or simply calculate as follows: Let $P , Q$ be two $k$-safe operators with coefficients $P_j$ and $Q_j$, respectively. Then, for the coefficient of $A:= P \o Q$ belonging to the multiindex $J$ we get similarly

$$ A_J = \sum_{j_1 + j_2 \geq J \geq j_2} \binom{j_1}{j_1 + j_2 - J} P_{j_1} \partial_{j_1 + j_2 - J} Q_{j_2}  \in H^{k-(s_1 + s_2) +\vert J \vert} = H^{k - {\rm deg} A + \vert J \vert} , $$

\v thus $A$ is $k$-safe, and clearly the map is continuous. \hfill \qed

\bigskip

\begin{Theorem}
\label{Fredholm}

Let $P$ be an elliptic $l$-safe operator of order $s$ with $l > n/2$. 
\begin{enumerate}
\item If the formal adjoint of $P$ is $(l-s)$-safe, then $P^{(m)}$ is Fredholm for all $m \leq l$, and its index equals the topological index. 
\item If $P$ is formally self-adjoint, then ${\rm ind} (P^{(l)}) =0$.
\end{enumerate}
\end{Theorem}

\v {\bf Proof.} Theorem \ref{CStarAlgebra} (Item 2) and the fact that $S_{k+1}^s \subset S_k^s$ imply that we can w.l.o.g. focus on the case $m=l$. To show that $P^{(l)}: H^l \rightarrow H^{l-s}$ is Fredholm, the main tool is the following elliptic estimate for nonsmooth coefficients proven in the Appendix:

\begin{Lemma}[Nonsmooth elliptic regularity]
\label{garding}
Let $k>n/2$, and let $P$ be an elliptic $k$-safe operator of order $s$, then there is $C>0$ such that for every $q \in \N$ with $q \leq k$ and every $u \in H^{q+s-1}$ with $Pu \in H^q$ we have $ u \in H^{q+s}$, and $\vert \vert u \vert \vert_{q+s} \leq C ( \vert \vert Pu \vert \vert_q + \vert \vert u \vert \vert_{q+s-1})  $.
\end{Lemma}

\v We also use the following quite classical closed range lemma whose proof we leave to the reader:

\begin{Lemma}
Let $X,Y,Z$ be Banach spaces, let $P: X \rightarrow Y, K: X \rightarrow Z$ be bounded linear operators, $K$ compact. If there is a constant $C>0$ such that for each $x \in X$ we have

$$ \vert \vert x \vert \vert_X \leq C \big( \vert \vert P x \vert \vert_Y + \vert \vert K x \vert \vert_Z  \big) , $$

then $P(X)$ is a closed subspace of $Y$, and $\ker (P)$ is finite-dimensional.
\end{Lemma}

The lemma applied to $P$ and $q = s-l $ and to the compact canonical inclusion of $H^{q+s}$ in $H^{q+s-1}$ yields that $P$ has closed image and finite-dimensional kernel.
For the corresponding statement on the cokernel, we use $\coker(P) := P(H^{l})^\perp = \ker((P'_{g, \omega})^{(l-s)})$ for the adjoint $(P'_{g, \omega})^{(l-s)}: H^{l-s} \rightarrow H^{l-2s}$ and the fact that for an elliptic operator $P$ its formal adjoint $P'_{g, \omega}$ is elliptic again, due to the relation $\s (P') = (-1)^s \s (P)^g$ between the principal symbols $\s(P)$ of $P$ and $\s(P')$ of $P'$. Thus we can apply Lemma \ref{garding} and the closed range lemma to the formal adjoint $P'$ to show the statement on the cokernel in $H^{l-s}$.


Using that $C^{\infty }$ is dense in any Sobolev space, we approximate every local coefficient $P_{i,t}$ of $P$ by sequences $P_{i,t} (n)$ of smooth coefficients converging in $H^{a(i,k)}$. The (by means of a smooth partition of unity) so defined differential operators $P(n) $ are elliptic for sufficiently large $n$, and they converge to $P$ in the operator norm from $H^l$ to $H^{l-s}$ due to Theorem \ref{CStarAlgebra}, first assertion. 
As $P^{(l)}$ is Fredholm, and as the index is locally constant in the open subset of Fredholm operators with the subspace topology, application of the classical index theorem for smooth coefficients completes the proof of the first assertion. 


For the second assertion, use that Fredholmness implies the last equality in  

$\dim \coker P^{(l)} = \dim \ker (P')^{(l-s)}  = \dim \ker P^{(l-s)}  = \dim \ker P^{(l)} .$ \hfill \qed

\bigskip

Let us consider examples of $k$-safe operators.

\begin{Theorem}
\label{Bochner}
Let $(\pi: E \rightarrow M, \g, \nabla)$ be a metric $H^k$ vector bundle with connection over $M$. Then $\nabla: \Gamma (\pi) \rightarrow \Gamma(\tau ' M) \x \Gamma(\pi)$ is $k$-safe. If moreover the smooth tangent bundle bundle $\tau M:TM \rightarrow M$ is equipped with a metric $H^{k-1}$ metric for $k>\frac{n}{2} +1$, then the Bochner Laplacian $\Delta^\pi$ on $\pi$ is $k$-safe. \end{Theorem}

\V {\bf Proof.} The operator $\nabla$ is $k$-safe by counting orders: The first coefficients are smooth (and constant) and the second are in $H^k$. Now, $ \tau ' M \x \pi $ is a metric $H^{k-1}$-vector bundle. Concretely, its connection $\tilde{\nabla}: \Gamma (\tau ' M \x \pi ) \rightarrow \Gamma ((\tau'M)^2 \x \pi)$ maps $b \x f$ to $\sum a_i \x (\nabla_{a_i^g } b \x f + b \x \nabla_{a_i^g} f)$, where $a_i$ is a frame for $ \tau' M$ and $ a_i^g $ denotes the $g$-equivalent vector field.
The regularities of the coefficients $P_2, P_1, P_0$ of $P:= \tilde{\nabla} \o \nabla$ are $p_2 = \infty, p_1 = k-1, p_0 = k-1$. As $tr^g = g^{-1} \o $, which is an $H^{k-1} $ endomorphism (by the cofactor formula for the inverse matrix and e.g. the analogue of Fa\`a di Bruno's formula for higher (weak) derivatives of the reciprocal of a function, always keeping in mind that $g$ is continuous and thus is uniformly positive by compactness), the corresponding regularities $q_2, q_1, q_0$ of $ Q:= \Delta^\pi := \tr^g ( \tilde{\nabla}  \o \nabla)$ are $q_2 = k -1, q_1 = k-1, q_0 = k-1$, in particular $q_j \geq j - 2 + \lfloor n/2 \rfloor +1$ and $q_j \geq k - 2$, thus $\Delta^\pi$ is $k$-safe due to Theorem \ref{CStarAlgebra}.  \hfill \qed

\begin{Theorem}
\label{pullback}
Let $(N, \overline{g})$ be a Riemannian manifold, $f \in H^k (M,N)$ an immersion. Then $(f^* \tau N, f^* \overline{g}, f^* \nabla)$ is a metric $H^k$ vector bundle with connection, whereas $(\tau M, h:=  f^* \overline{g}, \nabla^h)$ is a metric $H^{k-1}$ vector bundle with connection. Moreover, the second fundamental form of $f$ is an $H^{k-1}$ endomorphism. 
\end{Theorem}

\V {\bf Proof.} If we denote the bundle chart transitions of a bundle $\pi$ by $t_{ij}$, the bundle chart transitions of the pull-back bundle $f^* \pi$ are $t_{ij} \o f$. Let $S \in C^1 (f^* \tau N)$ and let $\{ x_i \} $ be local coordinates in a neighborhood $U$ of $q=f(y) \in f(M)$ with associated Christoffel symbols $\Gamma_{\gamma \a}^\b$ and the corresponding decomposition $S(y') = \sum_\a S_\a (y') \cdot  (\frac{\partial}{\partial x_\a} \o f) $ for $y'$ in the neighborhood $f^{-1} (U)$ of $y$. For $S \in H^l$, we have 

\begin{equation*}
(f^* \nabla )_X S = \big( \sum_\b X(S_\b) + \sum_{\gamma, \a, \b} \G_{\gamma \a}^\b \cdot (d_y f \cdot X)_\gamma \cdot S_\a \big) \cdot (\frac{\partial}{\partial x_\b} \o f),
\end{equation*} 

and as $\G_{\gamma \a}^\b \in C^\infty (U)$ and $d_y f \cdot X \in H^{k-1}$, the Christoffel symbols $\tilde{\Gamma}_{i \a}^\beta := \sum \Gamma_{\gamma \alpha}^\beta \cdot \frac{\partial f_\gamma}{X_i}$ in these charts are of regularity $H^{k-1}$, thus the bundle $f^*\tau N$ is a metric $H^k$ vector bundle with connection. As for the second fundamental form, $T (f(N)) $ is an $H^k$ submanifold of $TN$, and the orthogonal projection $p$ to its pullback in $f^* \tau N $ is an $H^k$ endomorphism. Thus $S:= p \o f^* \nabla$ is a $k$-safe operator of order $1$ and an endomorphism of regularity $H^{k-1}$ . \hfill \qed 

\bigskip

We write $g := f ^* \overline{g}$ on $f^* \tau N$, and for $m \in \N$ and $C \in \R$ we define $P_{C,m} := \Delta^g - C\cdot Q_m$ where $C$ is a constant and $Q_m$ is the $m$-th invariant polynomial of the second fundamental form of $f$, e.g.  $Q_0:=  \vert {\rm tr}_g S^f \vert^2 $, the mean curvature. Then Theorems \ref{Fredholm}, \ref{Bochner} and \ref{pullback} imply

\begin{Theorem}
\label{Application}
Let $(N, \overline{g})$ be a Riemannian manifold, $k>n/2 $  and $f \in H^k (M,N)$. 
If $l \in \N$, $2 \leq l$, then $  P_{C,m} \vert_{H^l}$ is a Fredholm operator from $H^l$ to $H^{l-2}$ for all $l \leq k$ and its index is zero. Moreover, $P_{C,m}$ is a smooth map from $Imm^k(M,N)$ to $ H^l(f^*\tau N) ' \x  H^{l-2}(f^* \tau N)$ for all $2 \leq l \leq k$. \hfill \qed
\end{Theorem}

\V In order to illustrate the subtleties of the use of the index theorem in the context of the geodesic equation on Hilbert manifolds, let us consider the remarkable article \cite{BHM} of Bauer, Harms and Michor on spaces of mappings and shape spaces. As noted by the authors themselves, the proof of their statement on Page 415 about local well-posedness of the geodesic equation for a metric defined by an operator $P$ is not entirely correct, as the assumptions 1 - 3 on $P: f \mapsto P_f$ on Page 414 concern exclusively smooth maps $f$, thus it is unclear whether each such $P$ necessarily extends continuously from the space of smooth immersions between two manifolds $M$ and $N$ to ${\rm Imm}^{k+2q} (M,N)$ taking at one such immersion $f$ a value $P_f$ in the bounded linear operators from $T_f {\rm Imm}^{k+2q} (M,N) = H^{k+2q} (f^* \tau N)$ to $ T_f {\rm Imm}^k (M,N)  = H^k (f^* \tau N)$ (consider e.g. $P_f:= (1 + Q_0) \cdot \Delta^{f^* \tau N}$). To make things easy, let us just include in the assumptions that {\em $P$ extends continuously to each Sobolev manifold ${\rm Imm}^{k+2q} (M,N)$, and that this extension is a $(k+2q)$-safe operator from each $T_f{\rm Imm}^{k+2q} (M,N)$ to $T_f {\rm Imm}^k (M,N)$}. As shown in Theorem \ref{Application}, this assumption is satisfied by the principal example $\1 + \Delta^q$ of that article. Let us call the local well-posedness of geodesics with this new assumption on $P$ the {\em modified BHM Theorem}. In the continuation of the proof, the authors refer to Theorem 26.2 of Shubin's book \cite{mS} to show vanishing index (implying then bijectivity) of any $P_f$. But the theorem in this reference applies to operators with smooth coefficients only, using explicitly that the operator maps smooth sections to smooth sections. Thus the theorem in Shubin's book cannot be applied in the present situation. Rather, one should use Theorem \ref{Application} for $l=d=k+2q$, completing the proof of the modified BHM Theorem. In a nutshell, the gap in the cited article can be closed by the additional assumption of safeness, which, due to Theorems \ref{CStarAlgebra} and \ref{Application}, holds true in the example $P:= \1 + \Delta^q$, inducing in fact a $C^\infty$ map 

$$P : \Imm^{k+2q} \ni f \mapsto P_f \in L(T_f \Imm^{k+2q} , T_f \Imm^k ).$$

\bigskip

Let us remark that the results obtained in this article could be extended in various ways:

\begin{enumerate}
\item Given sufficient regularity, using difference operators it is easy to see that even the stronger estimate $\vert u \vert_{l+s} < C (\vert Pu \vert_l + \vert u \vert_l)$ holds in Lemma \ref{garding}.
\item Everything could be done for nonempty boundary and Dirichlet or Neumann boundary conditions.
\item One can generalize the definition of $k$-safeness to pseudodifferential operators such that the statements above still hold.
\item Using the treatment of multiplicativity properties of negative Sobolev spaces as e.g. in \cite{BH} one could extend the $k$-safeness condition correspondingly. 
\item Finally, requiring appropriate decay conditions and using carefully modified Sobolev ring properties one could transfer the statement to the noncompact setting.
\end{enumerate}

\section{Appendix: Elliptic estimates for Sobolev coefficients}

Now we want to the prove Lemma \ref{garding}, which is to be performed in close analogy to the proof for the classical standard elliptic estimates in the smooth case as laid down for example in Taylor's book \cite{mT}, Th. 5.11.1., one of the most foundational results of elliptic theory, and it is quite likely that the case of Sobolev coefficients has been considered already, yet the author could not find a version of the lemma in the literature.  

First of all, pick a covering $U$ of $M$ by finitely many open chart neighborhoods $U_1, ... , U_n$ with charts $\kappa_1, ... , \kappa_n$. W.l.o.g. the chart neighborhoods $\kappa_i(U_i) $ are precompact domains in $\R^n$. Choose a smooth partition of unity $\psi_i$ subordinate to $U$. Then the operator $Q_i:= [P, \psi_i \cdot ]$ is of order $s-1$, and Theorem \ref{CStarAlgebra} implies that it is $k$-safe.

It is easy to show the estimate of the lemma for $p=0$: Ellipticity (entailing uniform ellipticity as $M$ is compact) implies that $\vert \vert \nabla^i u \vert \vert_{L^2} \leq C \vert \vert Pu \vert \vert_{L^2}$ if $\vert i \vert =s$ (cf. \cite{mT}), and application of the interpolation estimates yields $\vert \vert u \vert \vert_{H^s} \leq C \vert \vert Pu \vert \vert_{L^2} + \vert \vert u \vert \vert_{L^2}$ (cf. \cite{mT}). 

For $p>0$, as a first step we want to prove that we can {\em localize} the statement in the sense that it suffices to show it for ${\rm supp} (u) \in V$, where $V$ is an open precompact set in $\R^n$. To see this, suppose we know the estimate for that case, then (the $\psi_i$ being elements of a partition of unity and the $C_n$ constants independent of $u$) we get

\bean
\vert \vert u \vert \vert_{q+s} &\leq& C_1 \sum_i \vert \vert \psi_i u \vert \vert_{q+s} \leq C_2 \sum_i \big( \vert \vert P \psi_i u  \vert \vert_q + \vert \vert  \psi_i u \vert \vert_{q+s-1}  \big)  \\
&\leq&  C_3 \sum_i \big( \vert \vert \psi_i Pu \vert \vert_q + \vert  \vert Q_i  u \vert \vert_q + \vert \vert \psi_i u \vert \vert_{q+s-1} \big)   \leq C_4 (\vert \vert P u \vert \vert_q + \vert \vert u \vert \vert_{q+s-1}).
\eean

Thus, indeed, the statement for compactly supported operators in $\R^n$ implies the statement on compact manifolds. Consequently, in the following we suppress the $t$-indices in the coefficients $P_i$.

Our ansatz to prove the statement on $V \subset \R^n$ is to freeze the coefficients, like in Taylor's proof. But evaluating the coefficients on a lattice, as it is done in \cite{mT}, would not be a well-defined operation in our case, unless the coefficient $P_i$ in question is continuous, that is to say, for $\vert i \vert = s$. Instead of evaluating we will consider means of coefficients, i.e. associate for $z \in \Z^n$ and the point $\e z$ on the lattice $\e \Z^n$ the value $\hat{P}_{i, \e} (\e z) =  \fint_{B_\e (\e z)} P_i =  {\rm vol}(B_\e (\e z))^{-1} \cdot \int_{B_\e (\e z)} P_i$, and we define $P_{  \e} (y) := \sum_{\vert i \vert \leq s} \hat{P}_i(y) \partial_i $. Choose a partition of unity in $\R^n$ of the form $\{ \phi_z = \phi_0 (\cdot - z) \vert z \in \Z^n  \} $, for a smooth compactly supported function $\phi_0$ on $\R^n$. Such a partition can easily be constructed by pulling back a smooth partition of unity on the torus $\mathbb{T}^n$ subordinated to a trivializing atlas of the universal covering $\R^n \rightarrow \mathbb{T}^n$. Then define partitions of unity $\{ \phi_{\e z, \e}  \vert z \in \Z^n \} $ in $\R^n$ centered at the rescaled lattice $\e \Z^n$ by $\phi_{\e z , \e} (x) := \phi_0 (\frac{x- \e z}{\e}) = \phi_0 (\frac{x}{\e} - z) = \phi_z (x/ \e) $. In the end we will use $ \lim_{\e \rightarrow 0} \sum_{y \in \e \Z^n} \hat{P}_{i, \e} (y) \phi_{y, \e}  = P_i$ in $C^0$ for $\vert i \vert = s$, as those coefficients are continuous. The constant coefficient case is then treated by standard Fourier methods. Now $Pu = f$ implies

\bean
\phi_{\e z, \e} P_{\e} (\e z) u = \phi_{\e z , \e} f - R_{z, \e} u
\eean

where

\bea
\label{uniform}
R_{z, \e} (x) := \phi_{\e z, \e } \sum_{\vert i \vert \leq s} \big( P_i (x) - \hat{P}_{i, \e} (\e z) \big) \partial_i
\eea

Thus, with $Q_{ z, \e} := [P_{ \e} (\e z), \phi_{\e z, \e}]$, which is a smooth differential operator of order $s-1$, we get

\bea
\label{freezing}
P_{ \e } (\e z) (\phi_{\e z, \e} u) = \phi_{\e z, \e} f - R_{z, \e} u - Q_{   z, \e} u .
\eea

Now, for a differential operator $A$ on $\R^n$ let $\tilde{A}:= \mathcal{F}^{-1} \o A \o \mathcal{F}$, where $\mathcal{F}$ is the Fourier transform. Ellipticity of $P$ means that $\widetilde{P}_{\e}  (\e z) \vert_{\R^n \setminus B_K(0)} $ is bounded from zero for all $K >0$. Therefore homogeneity of $\widetilde{P}_{\e z, \e}$ implies that we can pick $\kappa \in C^\infty (\R^n, \R)$ of compact support and with $\kappa (B_K(0)) = \{ 1 \} $ and define a pseudodifferential operator $E_{z, \e }$ by

\bean
\widetilde{E_{ z, \e}} := (1 - \kappa) \widetilde{P_{ \e} (\e z)}^{-1} 
\eean

Then there is a $C>0$ independent of $\e$ such that $\vert \widetilde{E_{ z , \e}} (\xi) \vert \leq C (1 + \xi)^{-s}$, and thus , for all $\sigma \in \N$, $E_{ z , \e}: H^\s \rightarrow H^{\s + s}$ is continuous with operator norm bounded uniformly in $\e$:

\bean
\exists K_\s >0 \ \  \forall \e >0 \ \ \forall v \in H^\s : \vert \vert v \vert \vert_{\s + s } \leq K \cdot \vert \vert v \vert \vert_\s . 
\eean

Furthermore, we have

\bea
\label{almostinverse}
E_{z, \e} \o P_{ \e} (\e z) = \1 + \rho_{ z, \e} 
\eea 

for a pseudodifferential operator $\rho_{ z , \e}$ with $\widetilde{\rho_{ z , \e}} \in C_c^{\infty}$, thus $\rho_{z, \e}$ is infinitely smoothing, that is, maps any Sobolev space $H^\s$ to $C^\infty$, and its operator norms $H^\s \rightarrow H^\tau$ are bounded for all choices of $\s$ and $\tau$. Now we apply $E_{z, \e}$ to Equation (\ref{freezing}), insert Equation (\ref{almostinverse}) and take the sum over all $z \in \Z^n$, obtaining

\bean
\label{splitting}
u= \sum_{z \in \Z^n} E_{ z , \e} (\phi_{\e z , \e} f) + E_{z, \e }R_{z, \e} u + E_{\e z , \e}Q_{z, \e } u + \rho_{z, \e} (\phi_{\e z , \e} u)
\eean

Now it's all about estimating the four terms on the right-hand side. The last one is unproblematic as it contains an infinitely smoothing operator. For the first one, using that $\phi $ is smooth and compactly supported and counting orders, we obtain

\bean
\vert \vert \sum_{z \in \Z^n} E_{ z , \e}(\phi_{\e z, \e} f) \vert \vert_l \leq C(l , \e) \vert \vert f \vert \vert_{l-s} .
\eean

Taking into account that ${\rm deg}(Q_{z, \e}) < s-1$, the two last terms can be estimated as

\bean
\vert \vert  \sum_{z \in \Z^n} E_{ z , \e} Q_{z, \e} u + \rho_{ z, \e} (\phi_{\e z , \e} u) \vert \vert_l \leq C(l, \e) \vert \vert u \vert \vert_{l-1} .
\eean

Finally, looking at Equation (\ref{uniform}), using that the highest coefficients of $P$ are continuous (as $k>n/2$) and thus uniformly continuous as $M$ is compact, for all $l \leq k$ the second term can be estimated as

\bean
\vert \vert \sum_{z \in \Z^n} E_{ z , \e} R_{ z, \e} u \vert \vert_l &\leq& K_l \vert \vert R_{z, \e} u \vert \vert_{l-s}\\
 &=& K_l \cdot \vert \vert \sum_{z \in \Z^n} \Phi_{\e z, \e } \sum_{\vert i \vert \leq s } (P_i (x) - \hat{P}_{i, \e} (\e z ) ) \partial_i u \vert \vert_{l-s}\\
 &=& K_l \cdot \vert \vert \sum_{z \in \Z^n} \sum_{\vert k \vert \leq l-s} \partial_k ( \Phi_{\e z, \e } \sum_{\vert i \vert \leq s} (P_i (x) - \hat{P}_{i, \e} (\e z ) ) \partial_i u)  \vert \vert_0\\
 &\leq& K_l \cdot \vert \vert \sum _{z \in \Z^n} \sum_{\vert k \vert = l-s, \vert i \vert = s}  \Phi_{\e z, \e } \cdot (P_i (x) - \hat{P}_{i, \e} (\e z ) ) \partial_i \partial_k u \vert \vert_0 +  C(l, \e) \vert \vert u \vert \vert_{l-1}\\ 
&\leq& C(l) A(\e) \vert \vert u \vert \vert_l + C(l, \e) \vert \vert u \vert \vert_{l-1} ,
\eean

with a constant $C(l)$ independent of $\e$ and with $A(\e) = \vert \vert \sum_{\vert i \vert = s }\sum_{z \in \Z^n} \Phi_{\e z, \e } \cdot (P_i (x) - \hat{P}_{i, \e} (\e z) )   \vert \vert_{\infty}$. Now, as $\vert i \vert = s$ means considering the highest coefficients, which are continuous, we have ${\rm lim}_{\e \rightarrow 0} A(\e) = 0$. Thus we can choose $\e < (2 C(l))^{-1}$ and incorporate the second term into the left-hand side yielding a factor $1/2$, by which the desired form of the estimate is obtained. \hfill \qed





\end{document}